\magnification 1200
\advance\hoffset by 1,5truecm
\advance\hsize by 2,3 truecm
\def\makefootline{\baselineskip=52pt\line{\the\footline}}
\vsize= 23 true cm
\hsize= 14 true cm
\overfullrule=0mm
\font\sevensl=cmti7

\headline={\hfill\tenrm\folio\hfil}
\footline={\hfill}\pageno=1

\newcount\coefftaille \newdimen\taille
\newdimen\htstrut \newdimen\wdstrut
\newdimen\ts \newdimen\tss

\long\def\demdrien#1{{\parindent=0pt\messages{debut de preuve}\smallbreak
     \advance\margeg by 2truecm \leftskip=\margeg  plus 0pt
     {\everypar{\leftskip =\margeg  plus 0pt}
              \everydisplay{\displaywidth=\hsize
              \advance\displaywidth  by -1truecm
              \displayindent= 1truecm}
     {\bf D\'emonstration } -- \enspace #1
      \hfill}\bigbreak}\messages{fin de preuve}}

\long\def\demdriensanstext#1{{\parindent=0pt\messages{debut de preuve}\smallbreak
     \advance\margeg by 2truecm \leftskip=\margeg  plus 0pt
     {\everypar{\leftskip =\margeg  plus 0pt}
              \everydisplay{\displaywidth=\hsize
              \advance\displaywidth  by -1truecm
              \displayindent= 1truecm}
     {\bf D\'emonstration } -- \enspace #1
      \hfill}\bigbreak}\messages{fin de preuve}}

\def\fixetaille#1{\coefftaille=#1
\htstrut=8.5pt \multiply \htstrut by \coefftaille \divide \htstrut by 1000
 \wdstrut=3.5pt \multiply \wdstrut by \coefftaille \divide \wdstrut by 1000
\taille=10pt  \multiply \taille by \coefftaille \divide \taille by 1000
\ts=\taille \multiply \ts by 70 \divide \ts by 100
 \tss=\taille \multiply \tss by 50 \divide \tss by 100
\font\tenrmp=cmr10 at \taille
\font\sevenrmp=cmr7 at \ts
\font\fivermp=cmr5 at \tss
\font\tenip=cmmi10 at \taille
\font\sevenip=cmmi7 at \ts
\font\fiveip=cmmi5 at \tss
\font\tensyp=cmsy10 at \taille
\font\sevensyp=cmsy7 at \ts
\font\fivesyp=cmsy5 at \tss
\font\tenexp=cmex10 at \taille
\font\tenitp=cmti10 at \taille
\font\tenbfp=cmbx10 at \taille
\font\tenslp=cmsl10 at \taille}

\def\fspeciale{\textfont0=\tenrmp%
\scriptfont0=\sevenrmp%
\scriptscriptfont0=\fivermp%
\textfont1=\tenip%
\scriptfont1=\sevenip%
\scriptscriptfont1=\fiveip%
\textfont2=\tensyp%
\scriptfont2=\sevensyp%
\scriptscriptfont2=\fivesyp%
\textfont3=\tenexp%
\scriptfont3=\tenexp%
\scriptscriptfont3=\tenexp%
\textfont\itfam=\tenitp%
\textfont\bffam=\tenbfp%
\textfont\slfam=\tenbfp%
\def\it{\fam\itfam\tenitp}%
\def\bf{\fam\bffam\tenbfp}%
\def\rm{\fam0\tenrmp}%
\def\sl{\fam\slfam\tenslp}%
\normalbaselineskip=12pt%
\multiply \normalbaselineskip by \coefftaille%
\divide \normalbaselineskip by 1000%
\normalbaselines%
\abovedisplayskip=10pt plus 2pt minus 7pt%
\multiply \abovedisplayskip by \coefftaille%
\divide \abovedisplayskip by 1000%
\belowdisplayskip=7pt plus 3pt minus 4pt%
\multiply \belowdisplayskip by \coefftaille%
\divide \belowdisplayskip by 1000%
\setbox\strutbox=\hbox{\vrule height\htstrut depth\wdstrut width 0pt}%
\rm}

\def\vmid#1{\mid\!#1\!\mid}

\def\fle{\rightarrow}

\null\vskip-1cm
\font\sevensl=cmti7
\font\sc=cmcsc10

\newdimen\emm 
\def\pmb#1{\emm=0.03em\leavevmode\setbox0=\hbox{#1}
\kern0.901\emm\raise0.434\emm\copy0\kern-\wd0
\kern-0.678\emm\raise0.975\emm\copy0\kern-\wd0
\kern-0.846\emm\raise0.782\emm\copy0\kern-\wd0
\kern-0.377\emm\raise-0.000\emm\copy0\kern-\wd0
\kern0.377\emm\raise-0.782\emm\copy0\kern-\wd0
\kern0.846\emm\raise-0.975\emm\copy0\kern-\wd0
\kern0.678\emm\raise-0.434\emm\copy0\kern-\wd0
\kern\wd0\kern-0.901\emm}

\font\tendb=msbm10
\font\sevendb=msbm7

\newfam\dbfam
\textfont\dbfam=\tendb\scriptfont\dbfam=\sevendb\scriptscriptfont\dbfam=\sevendb
\def\db{\fam\dbfam\tendb}

\def\C{{\db C }}

\def\N{{\db N }}

\def\R{{\db R }}

\def\Z{{\db Z }}


\font\gothique=eufm10
\def\got#1{{\gothique #1}}

\newdimen\margeg \margeg=0pt
\def\bb#1&#2&#3&#4&#5&{\par{\parindent=0pt
    \advance\margeg by 1.1truecm\leftskip=\margeg
    {\everypar{\leftskip=\margeg}\smallbreak\noindent
    \hbox to 0pt{\hss\bf [#1]~~}{\bf #2 - }#3~; {\it #4.}\par\medskip
    #5 }
\medskip}}

\newdimen\margeg \margeg=0pt
\def\bbaa#1&#2&#3&#4&#5&{\par{\parindent=0pt
    \advance\margeg by 1.1truecm\leftskip=\margeg
    {\everypar{\leftskip=\margeg}\smallbreak\noindent
    \hbox to 0pt{\hss [#1]~~}{\pmb{\sc #2} - }#3~; {\it #4.}\par\medskip
    #5 }
\medskip}}

\newdimen\margeg \margeg=0pt
\def\bba#1&#2&#3&#4&#5&{\par{\parindent=0pt
    \advance\margeg by 1.1truecm\leftskip=\margeg
    {\everypar{\leftskip=\margeg}\smallbreak\noindent
    \hbox to 0pt{\hss [#1]~~}{{\sc #2} - }#3~; {\it #4.}\par\medskip
    #5 }
\medskip}}

\def\messages#1{\immediate\write16{#1}}

\def\findem{\vrule height0pt width4pt depth4pt}

\long\def\demA#1{{\parindent=0pt\messages{debut de preuve}\smallbreak
     \advance\margeg by 2truecm \leftskip=\margeg  plus 0pt
     {\everypar{\leftskip =\margeg  plus 0pt}
              \everydisplay{\displaywidth=\hsize
              \advance\displaywidth  by -1truecm
              \displayindent= 1truecm}
     {\bf Proof } -- \enspace #1
      \hfill\findem}\bigbreak}\messages{fin de preuve}}

\def\resp{\mathop{\rm resp}\nolimits}
\def\resp.{\mathop{\rm resp.}\nolimits}

\null\vskip 1,5cm
\centerline{\bf Fonctions holomorphes Cliffordiennes}
\bigskip
\centerline{{\bf Guy Laville } et {\bf Ivan Ramadanoff}}

\bigskip\medskip
{\leftskip=0cm\rightskip=10pt
{\parindent=1cm\narrower\fixetaille{700}{\fspeciale {\bf R\'esum\'e.-} \  {\sevenrm Soit} \
$\R_{0,2m+1}$ \ {\sevenrm l'alg\`ebre de Clifford de $\R^{2m+1}$ muni d'une forme quadratique
de signature n\'egative,}
$D =
\displaystyle\sum_{i=0}^{2m+1} \ e_i \ \displaystyle{\partial\over \partial \
x_i}$, \ $\Delta$ \ {\sevenrm le Laplacien ordinaire. Les fonctions holomorphes
Cliffordiennes} $f$  {\sevenrm sont les fonctions satisfaisant \`a} \ $D\Delta^m f = 0$.
{\sevenrm Nous \'etudions les solutions polynomiales et singuli\`eres, les
repr\'esentations int\'egrales et leurs cons\'equences et enfin le fondement de
la th\'eorie des fonctions elliptiques Cliffordiennes.}}
\par}

\bigskip\bigskip
\centerline{{\bf Cliffordian holomorphic functions}}
\medskip
{\parindent=1cm\narrower\fixetaille{700}{{\fspeciale {\bf Abstract.-} \  {\sevensl Let} \
$\R_{0,2m+1}$
\ {\sevensl be the Clifford algebra of} \ $\R^{2m+1}$  {\sevensl with a quadratic form  of
negative signature,} \ $D = \displaystyle\sum_{i=0}^{2m+1} \ e_i \
\displaystyle{\partial\over \partial x_i}$, \ $\Delta$  {\sevensl the ordinary Laplacian.
The holomorphic Cliffordian functions are solutions of} \ $D\Delta^m f = 0$.  {\sevensl We
study the polynomial and singular solutions, re\-presentation integral formulas
and the foundation of the Cliffordian elliptic function theory.}}
}
\par}

\bigskip\bigskip
\bigskip\bigskip
\parindent=0,5cm
{\sc 1. Introduction.} -- La th\'eorie des fonctions d'une variable complexe 
(dimension r\'eelle 2) a \'et\'e d\'evelopp\'ee en dimensions sup\'erieures
dans deux directions~: les fonctions holomorphes de plusieurs variables
complexes (utilisation du corps $\C$, dimension r\'eelle $2n$) et la th\'eorie
des fonctions monog\`enes  (\'equation de Dirac,  utilisation des alg\`ebres de
Clifford, dimension r\'eelle $n$). Aucune de ces deux th\'eories n'est
satisfaisante si l'on veut avoir certaines parties tr\`es f\'econdes de la
premi\`ere th\'eorie~: fonctions elliptiques, fonction th\'eta, etc ...  \ De
fa\c con analogue la th\'eorie des fonctions de ``plusieurs variables
Cliffordiennes"  pr\'esente certaines difficult\'es dues au manque de solutions
des op\'erateurs  (voir [2], [3], [4]).

\medskip
Ceci nous conduit de fa\c con naturelle \`a la th\'eorie pr\'esent\'ee.

\bigskip\bigskip
{\sc 2. Notations et d\'efinitions.} --  Soit \ $\R_{0,2m+1}$ \ l'alg\`ebre de
Clifford de l'espace vectoriel $V$ de dimension r\'eelle $2m+1$  muni d'une
forme quadratique de signature n\'egative.  Soit $S$ l'ensemble des scalaires
de \ $\R_{0,2m+1}$, \ $S$ peut \^etre identifi\'e avec \ $\R$. Soit \ $\{ e_i\}$
\ $i = 1,\ldots , 2m+1$ \ une base orthonormale de $V$ et posons $e_0 = 1$. 
D\'efinissons l'op\'erateur (dit de Cauchy, Fueter, Dirac, voir [1])
$$D = \sum_{i=0}^{2m+1} \  \ e_i \ {\partial\over \partial x_i} .$$

\bigskip\bigskip
{\sc   D\'efinition 2.1.} -- \ {\sl Soit \ $\Omega$ \ un ouvert de $S\oplus V$.
Une fonction \ $f : \Omega \fle \R_{0,2m+1}$ \  sera dite holomorphe
Cliffordienne \`a gauche quand
$$D \Delta^m f = 0.$$

$\Delta^m$ \'etant le laplacien ordinaire it\'er\'e $m$ fois.
}

\bigskip\bigskip
{\sc Remarque 1} -- On pourrait ne consid\'erer que les fonctions $f :
\Omega\fle S \oplus V$  satisfaisant \`a l'\'equation ci-dessus.  Ces
fonctions peuvent engendrer les pr\'ec\'edentes par combinaisons lin\'eaires
\`a droite.

\bigskip\bigskip
{\sc Remarque 2} --  Soit \ $x = \displaystyle\sum_{i=0}^{2m+1} x_i e_i$,
$x_i\in\R$  alors $f$ est holomorphe Cliffordienne si et seulement si $\Delta^m
f(x) = 0$ et $\Delta^{m+1} \bigl( x f(x)\bigr) = 0$.

\bigskip\bigskip\bigskip
{\sc 3.  Solutions \'el\'ementaires} --  Posons :
$$\alpha = (\alpha_0,\ldots ,
\alpha_{2m+1}) \ \hbox{ avec } \ \alpha_i\in\N  \ \hbox{ et }  \ \vmid{\alpha}
\ = \sum_{i=0}^{2m+1} \alpha_i.$$  
Consid\'erons l'ensemble form\'e des \'el\'ements $\alpha_0$  fois $e_0$, \
$\alpha_1$ fois $e_1,\ldots ,
\alpha_{2m+1}$ fois \ $e_{2m+1}$~; \ cet ensemble sera not\'e \ $\{ e_\nu\}$. 
Posons~:
$$P_\alpha (x) = \displaystyle{ \ 1 \ \over \ \vmid{\alpha} ! \ } \
\sum_{\hbox{\got G}} \ \prod_{\nu = 1}^{\vmid{\alpha} -1} \ \bigl( e_{\sigma
(\nu )} x\bigr) \ e_{\sigma (\vmid{\alpha})}$$
la somme \'etant \'etendue sur tous les \'elements $\sigma$ du groupe des
permutations \got G .

\bigskip\bigskip
{\bf Th\'eor\`eme 3.1.} -- {\it $P_\alpha (x)$  est un polyn\^ome de degr\'e \
$\vmid{\alpha} - 1$ en $x$, holomorphe Cliffordien et tout polyn\^ome
holomorphe Cliffordien est combinaison lin\'eaire (\`a coefficients \`a
droite) de tels polyn\^omes.
} 

\medskip
Soit $\beta$ un multiindice du m\^eme type que $\alpha$.  Posons :
$$S_\beta (x) = { \ 1 \ \over \ \vmid{\beta} ! \ } \ \sum_{\hbox{\got G}} \
\prod_{\nu = 1}^{\vmid{\beta} - 1} \ \bigl( x^{-1} e_{\sigma (\nu )}\bigr) \
x^{-1}.$$

\bigskip\bigskip
{\bf Th\'eor\`eme 3.2.} -- {\it $S_\beta (x)$  est une fonction holomorphe
Cliffordienne d\'efinie sur  $(S \oplus V) \setminus \{ 0 \}$.
}

\bigskip
\demdriensanstext{Les d\'emonstrations de ces deux th\'eor\`emes se font soit par
un calcul direct, soit par l'\'etablissement du lien entre la structure
alg\'ebrique et la d\'erivation. }

\bigskip\bigskip
{\bf Th\'eor\`eme 3.3 (fraction rationnelle)} -- {\it Soit $P$ un entier, \ $\{ a_p\}$, $p =
1,\ldots , P$, \ $\{ b_q\}$, $q = 1,\ldots , P+1$ \  avec $a_p$ et $b_q$  \'el\'ements de
$S\oplus V$.  Alors il existe une fraction rationnelle, holomorphe Cliffordienne hors de ses
singularit\'es ayant les
$\{ a_p\}$ parmi ses z\'eros et les $\{ b_q\}$  parmi ses p\^oles.
}

\bigskip\bigskip\bigskip
{\sc 4.  Repr\'esentation int\'egrale et formule de Taylor} -- 
\medskip
Posons \ $N(x) = (-1)^m \  \displaystyle{m+1\over 2^{2m+1} m! \ \pi^{m+1}} \ x^{-1}$.

\bigskip\bigskip
{\bf Th\'eor\`eme 4.1.} -- {\it Soit $f : \Omega \fle \R_{0, 2m+1}$ \ une fonction holomorphe
Cliffordienne et $\Gamma$ un ouvert born\'e \`a bord r\'egulier, $\overline\Gamma
\subset\Omega$. Alors :
$$\eqalign{f(x) &=\int_{\partial\Gamma} \ \Delta^m N(y-x) f(y) d\sigma (y)\cr
       &-\sum_{k=1}^m \ \int_{\partial\Gamma} \ {\partial\over \partial n} \Delta^{m-k}
            N(y-x) D\Delta^{k-1} f(y) d\sigma (y)\cr
       &+\sum_{k=1}^m \ \int_{\partial\Gamma} \Delta^{m-k} N(y-x) \  {\partial\over \partial
          n} \ D\Delta^{k-1} f(y) d\sigma (y).\cr}$$
}

\demdriensanstext{La d\'emonstration se fait par les m\'ethodes classiques de th\'eorie du
potentiel et on en d\'eduit le th\'eor\`eme suivant~:
}

\bigskip\bigskip
{\bf Th\'eor\`eme (formule de Taylor)} --  {\it Soit $f$ une fonction holomorphe Cliffordienne
au voisinage d'un point $a\in S \oplus V$.  Alors la s\'erie suivante est uniform\'ement
convergente dans un voisinage de $a$ et on a l'\'egalit\'e suivante dans ce voisinage
$$f(x) = \sum_\alpha \ P_\alpha (x-a) c_\alpha$$
la sommation \'etant faite pour tous les multiindices $\alpha$ et $c_\alpha \in \R_{0,2m+1}$.}

\bigskip\bigskip\bigskip
{\sc 5. D\'eveloppement de Laurent} --
\medskip
{\bf Th\'eor\`eme 5.1.} -- Soit $B$ une boule de centre $a$ et de rayon $R$ dans $S\oplus V$
et $f$ une fonction holomorphe Cliffordienne dans $B \setminus \{ a\}$.  Alors pour tout 
$x\in B \setminus \{ a\}$  on a~:
$$f(x) = \sum_a P_\alpha (x) c_\alpha + \sum_\beta S_\beta (x) d_\beta$$
o\`u \ $\alpha$ \ et \ $\beta$ \ sont de multiindices d\'efinis dans le paragraphe 3.

\bigskip
\demdriensanstext{La d\'emonstration de ce th\'eor\`eme est une application directe du
th\'eor\`eme de repr\'esentation int\'egrale.}

\bigskip\bigskip\bigskip
{\sc 6. Fonctions elliptiques Cliffordiennes}
\medskip
{\bf D\'efinition.} -- {\sl On appelle fonction elliptique Cliffordienne une fonction
holomorphe Cliffordienne $2m+2$ p\'eriodique d\'efinie sur $(S\oplus V) \setminus E$  o\`u $E$
est un ensemble de points tel que $E\cap K$  soit fini pour tout compact $K$.
}

\bigskip
Les th\'eor\`emes de la th\'eorie classique des fonctions elliptiques se g\'en\'eralisent
quand ceux-ci ne font intervenir implicitement ou explicitement que la structure vectorielle
de l'ensemble de ces fonctions.

\medskip
Construisons l'analogue de la fonction $\zeta$ de Weierstrass~: soient $\omega_0,\ldots ,
\omega_{2m+1}$, \ $2m+2$ \'el\'ements de $S\oplus V$, \ $\R$-lin\'eairement ind\'ependants.
Posons~:
$$\Omega_K = 2 \ \sum_{j=0}^{2m+1} k_j \omega_j \quad\hbox{ avec }\quad K = (k_0,\ldots ,
k_{2m+1}) \in \Z^{2m+2}.$$

\bigskip\bigskip
{\bf Th\'eor\`eme 6.1.} --   {\sl  La fonction
$$\zeta_{2m+2}(x) = x^{-1} + \sum_{K\not= 0} \ \Bigl( (x-\Omega_K)^{-1} + \sum_{p=0}^{2m+1}
(\Omega_K^{-1} x)^p \ \Omega_K^{-1}\Bigr)$$
est bien d\'efinie par une s\'erie uniform\'ement convergente sur tout compact ne contenant
pas les singularit\'es,  est holomorphe Cliffordienne hors de ses singularit\'es. Ses
d\'eriv\'ees d'ordre sup\'erieur ou \'egal \`a $2m+1$  sont des fonctions elliptiques
Cliffordiennes.}

\bigskip
Comme dans le cas classique, cette fonction peut \^etre mise comme fondement de la th\'eorie
des fonctions elliptiques Cliffordiennes  (voir [5]).

\bigskip\bigskip
{\bf R\'ef\'erences bibliographiques}
\medskip
\bb 1&F. BRACKS, R. DELANGHE, F. SOMMEN&Clifford analysis&Pitman (1982)& &

\bb 2&G. LAVILLE&Une famille de solutions de l'\'equation de Dirac avec champ
\'electroma\-gn\'etique quelconque&CRAS, t. 296, pp. 1029-1032 (1983)& &

\bb 3&G. LAVILLE&Sur l'\'equation de Dirac avec champ \'electromagn\'etique
quelconque&Lecture notes n$^\circ$1165, pp. 130-149 (1985)& &

\bb4&V.P. PALAMODOV&On ``holomorphic" functions of several quatermonic va\-riables  C.A. Aytama
(ed.)&Linear topological spaces and complex analysis II Ankara (1995), pp 67-77& &

\bb5&J. TANNERY, J. MOLK&El\'ements de la th\'eorie des fonctions
elliptiques&Gauthier-Villars (1893-1902)& &

\bigskip\bigskip

{\parindent=8cm
\item{UPRES-A 6081} D\'epartement de Math\'ematiques
\item{}Universit\'e de Caen
\item{}14032 CAEN Cedex France

\medskip
\item{}glaville@math.unicaen.fr
\item{}rama@math.unicaen.fr
\par}
}

\end